\documentclass[a4paper,12pt]{article}

\advance\textheight2.9cm
\advance\topmargin-2.0cm
\advance\textwidth2.6cm
\advance\evensidemargin-1.6cm

\usepackage{amsfonts}
\usepackage{amssymb}
\usepackage{amsmath}
\usepackage{mathtools}
\usepackage{bm}

\usepackage{array}
\usepackage{multirow}
\usepackage{graphicx}
\usepackage{longtable}
\usepackage{booktabs}
\usepackage{hhline}
\usepackage{caption}
\usepackage{colortbl}
\usepackage[table]{xcolor}
\usepackage{adjustbox}
\usepackage{lscape}

\usepackage{algorithm}
\usepackage{algpseudocode}

\usepackage[utf8]{inputenc}
\usepackage[english]{babel}
\usepackage{lmodern}
\usepackage{textcomp}
\usepackage{csquotes}
\usepackage{url}
\usepackage{cite}
\usepackage{ulem}
\usepackage{cancel}
\usepackage{multicol}
\usepackage{blindtext}
\usepackage{comment}
\usepackage{chngcntr}

\usepackage{tikz}
\usetikzlibrary{matrix,shapes,arrows,positioning,chains}

\usepackage{fancyhdr}
\usepackage{hyperref}
\usepackage{cleveref}

\fancypagestyle{plain}{%
  \fancyhf{}%
  \fancyfoot[C]{\bfseries \thepage}%
}

\usepackage{theorem}

\newtheorem{rem}{Remark}

\newcolumntype{?}{!{\vrule width 1pt}}
\newcolumntype{R}[2]{%
  >{\adjustbox{angle=#1,lap=\width-(#2)}\bgroup}%
  l%
  <{\egroup}%
}

\makeatletter
\newcounter{subsubsubsection}[subsubsection]
\def\subsubsubsectionmark#1{}

\def\subsubsubsection{%
  \@startsection{subsubsubsection}{4}{\z@}%
  {-3.25ex plus -1ex minus -.2ex}%
  {1.5ex plus .2ex}%
  {\normalsize\bfseries}%
}
\def\l@subsubsubsection{\@dottedtocline{4}{4.8em}{4.2em}}

\makeatother

\definecolor{ao(english)}{rgb}{0.0, 0.5, 0.0}

\def\Box{{\hbox{\raisebox{0.0em}{\rlap{$\sqcap$}}\kern0em%
            \raisebox{-0.0em}{$\sqcup$}}} }

\newenvironment{proof}{{\it Proof. }}%
{\nopagebreak\hspace*{0.5cm}\hfill$\Box$\vspace{0.5cm}}

\def\bepf{\begin{proof}}
\def\epf{\end{proof}}

\def\fct#1{\mathop{\rm #1}}

\def\opt{\fct{opt}}

\def\rec{\fct{rec}}

\def\quadc{(2)}
\def\cubic{(3)}

\def\bfi#1{\textbf{#1}}

\def\eeq{\end{equation}}
\def\lbeq#1{\begin{equation}\label{#1}}
\def\bary{\begin{array}}
\def\eary{\end{array}}

\def\D{\displaystyle}
\def\ol{\overline}

\def\wt{\widetilde}

\def\eps{\varepsilon}

\def\p{\mathcal{P}}

\numberwithin{equation}{section}

\newcommand{\bd}{\mathbf{d}}
\newcommand{\bx}{\mathbf{x}}
\newcommand{\btau}{\boldsymbol{\tau}}
\newcommand{\bz}{\mathbf{z}}

\newcommand{\bP}{\mathbf{P}}

\newcommand{\bg}{\mathbf{g}}
\newcommand{\by}{\mathbf{y}}
\newcommand{\bs}{\mathbf{s}}

\newcommand{\bu}{\mathbf{u}}

\newcommand{\bp}{\mathbf{p}}
\newcommand{\bb}{\mathbf{b}}

\begin{document}

\begin{center}

{\Large \bf Diagonal Hessian Approximation Based on Conjugacy Condition for Noisy Derivative-Free Optimization Problems in High Dimensions }

\vspace{0.5cm}

{\large \bf Morteza Kimiaei}
\centerline{\sl Fakult\"at f\"ur Mathematik, Universit\"at Wien}
\centerline{\sl Oskar--Morgenstern--Platz 1, A--1090 Wien, Austria}
\centerline{\sl email: morteza.kimiaei@univie.ac.at}
\centerline{\sl \url{http://www.mat.univie.ac.at/~kimiaei}}

\vspace{0.5cm}

{\large \bf Saman Babaie--Kafaki}
\centerline{\sl Faculty of Engineering, Free University of Bozen--Bolzano}
\centerline{\sl NOI Techpark, Via Bruno Buozzi 1, 39100 Bolzano (BZ), Italy}
\centerline{\sl email: saman.babaiekafaki@unibz.it}

\vspace{0.5cm}

\end{center}

\bfi{Abstract.} We consider large-scale noisy derivative-free optimization (DFO) problems in which only function values are available and gradient or subgradient information cannot be reliably estimated. Matrix-adaptation evolution strategies ({\tt MAES}) and their limited-memory variants are among the most robust DFO methods under noise; however, their performance may deteriorate when the noise level is large. In such regimes, sorting and selection may misidentify informative sampled points, making the recombination step less reliable and weakening the scaling information used by affine or matrix-adaptation mechanisms. This can substantially reduce the efficiency of {\tt MAES}-type methods, especially in high-dimensional settings.

To address this limitation, we propose a DFO method that replaces the full affine-scaling matrix with a diagonal approximation constructed from conjugacy-type conditions. The proposed mechanism does not attempt to estimate gradients, subgradients, or interpolation models, nor does it learn dense covariance information from noisy rankings. Instead, it uses consecutive normalized recombination displacements in a conservative diagonal update, thereby limiting the influence of unreliable selection information while preserving the derivative-free structure of the underlying evolutionary framework. As a result, the method is computationally cheaper than full matrix-adaptation schemes and limited-memory affine-scaling variants, while providing a stable scaling mechanism in noisy environments. Numerical experiments on noisy benchmark problems show that the proposed method is competitive with, and often more efficient than, {\tt MAES}-type baselines, particularly when the noise level is large and ranking-based selection becomes unreliable.


\vspace{0.2cm} {\em 2000 AMS Subject Classification: 90C53, 65K05.}

\hfill \today


\section{Introduction}

Effectively estimating first-order or second-order information in \textbf{Nonlinear Programming} (NLP) is of fundamental importance from both modeling and algorithmic perspectives. Nowadays, the emergence of high-dimensional NLP models in data mining and machine learning, which are often highly complicated, highlights the need for limited-memory algorithms that benefit from approximate versions of the original model without explicitly containing gradient or Hessian information. Hence, in recent decades, \textbf{Derivative-Free Optimization} (DFO) has attracted significant attention as an effective class of NLP methodologies.

Among various NLP models, unconstrained optimization is widely recognized as one of the most fundamental and frequently studied cases \cite{Andreibook}, mainly due to its extensive applications in many decision-making processes \cite{SBKRairo}. The general form of such a problem is given by
\begin{equation}\label{UP}
\min_{\bx \in \mathbb{R}^n} f(\bx),
\end{equation}
where $f$ denotes the objective (or cost) function, which may be nonsmooth, noisy, or available only through inexact evaluations in practical applications. To solve problem \eqref{UP}, most existing algorithms are designed within either the line search or the trust region frameworks \cite{Andreibook}. Typically, these numerical schemes are iterative procedures, and can be represented as
\begin{equation*}
\bx_0 \in \mathbb{R}^n,\qquad \bx_{k+1} := \bx_k + \bs_k,\qquad k \geq 0.
\end{equation*}
In line search approaches, the step is defined by $\bs_k := t_k \bd_k$, where $t_k > 0$ denotes an (approximate) step size and $\bd_k \in \mathbb{R}^n$ is a search direction. In gradient-based methods, $\bd_k$ is usually required to be a descent direction (cf. \cite{Andreibook}), i.e.,
\[
    \bd_k^\top \bg_k<0,
    \qquad
    \bg_k:=\nabla f(\bx_k),
\]
where $\bg_k$ denotes the exact gradient  at $\bx_k$. Meanwhile, in trust region methods, the step $\bs_k$ is obtained by
approximately solving a constrained subproblem based on a local model
of the objective function \cite{Andreibook}.


In the DFO setting, however, $\bg_k$ is unavailable, and therefore the descent property cannot usually be verified directly. Consequently, $\bd_k$ must be generated without explicit derivative information. Common derivative-free choices include coordinate directions, polling directions, random directions, or directions obtained from interpolation or model-based procedures. Coordinate and polling directions are simple and structured, whereas random directions are often attractive in high-dimensional or noisy settings because they explore the space without constructing an explicit gradient approximation. Thus, in DFO, $\bd_k$ should be interpreted more generally as a generated search direction or trial direction, rather than as a direction whose descent property is known a priori; for more details see \cite{AudH,ConSV}.

Existing studies indicate that the efficiency of a line search algorithm largely depends on the choice of the search direction. In gradient-based line search methods, a typical search direction is given by $\bd_k := -\bP_k^{-1} \bg_k$, where $\bP_k$ is usually an $n \times n$ matrix (often as an approximation of the Hessian). In DFO, however, since $\bg_k$ is unavailable, the underlying algorithm therefore relies on function-value information together with generated steps, directions, or sampled candidate points, without requiring explicit gradient or Hessian information. Consequently, $\bP_k$ should be interpreted as a derivative-free curvature surrogate or scaling matrix. Meanwhile, for large-scale problems commonly encountered in modern applications, the use of dense matrices $\bP_k$ is computationally expensive and often impractical. As a result, practical implementations typically rely on implicit second-order information or limited-memory approximations of $\bP_k$, although such representations generally come at the cost of reduced accuracy \cite{SBKRairo}.

Lack of derivative information becomes even more pronounced in noisy DFO problems, where the available objective function values are contaminated by deterministic or stochastic perturbations. In such settings, even when finite-difference or interpolation-based gradient approximations are formally available, their accuracy may deteriorate severely, depending on the noise level, the sampling radius, and the geometry of the sampled points. Indeed, small perturbations in the sampled function values can be amplified by differencing procedures, leading to unreliable gradient surrogates and, consequently, unstable search directions, e.g., see \cite{MADFO,VRDFON}. Increasing the sampling radius may reduce the relative effect of noise, but it also decreases the locality of the approximation; conversely, decreasing the sampling radius improves locality but may amplify the noise. This trade-off makes accurate gradient reconstruction particularly challenging in high-dimensional noisy problems.

Several classes of DFO solvers have been developed for unconstrained noisy DFO problems. Since the present work does not aim to survey these methods, we only recall some representative references. These include derivative-free line-search-based solvers \cite[Section 2.3.4]{LarMW}, derivative-free trust-region-based solvers \cite[Section 2.4]{LarMW}, direct-search solvers \cite[Section 2.1]{LarMW}, and matrix adaptation evolution strategies, such as the solvers proposed by Auger and Hansen \cite{AugH}, Loshchilov et al. \cite{Loshchilov2019}, Beyer \cite{Beyer2020}, Beyer and Sendhoff \cite{Beyer2017}, and Kimiaei and Neumaier \cite{MADFO,MATRS}. Historical and methodological overviews of DFO can be found in the books by Audet and Hare \cite{AudH} and Conn et al. \cite{ConSV}. For comparative studies of the behavior of such solvers in the noiseless case, see Rios and Sahinidis \cite{RioS}, and Kimiaei and Neumaier \cite{KimN}; for the noisy case, see Kimiaei \cite{VRDFON}, and Kimiaei and Neumaier \cite{MADFO}. Further useful references on noisy DFO include Berahas et al. \cite{BerBN}, Elster and Neumaier \cite{ElsN0}, Gratton et al. \cite{GraRVZ,Gratton2017,GraTT}, Huyer and Neumaier \cite{HuyN}, Lucidi and Sciandrone \cite{LucS}, Mor\'e and Wild \cite{Mor2012}, Powell \cite{Pow02,Pow08}, Shi et al. \cite{shi2021adaptive}, and Wild et al. \cite{Wild08}.

To the best of our knowledge, in the noisy DFO setting, there remains a need for scaling and curvature-surrogate mechanisms that do not rely directly on finite-difference gradient estimates, explicit interpolation models, or gradient displacement vectors. This is particularly important in large-scale noisy problems, where constructing reliable derivative approximations may be computationally expensive and highly sensitive to perturbations in the sampled function values. Instead, one may exploit information already generated by the optimization process itself, such as accepted steps, successful mutation directions, or, more generally, consecutive generated directions. Although such quantities are still influenced by noisy function evaluations through the underlying selection or acceptance mechanism, they are not obtained from direct differencing of noisy objective values and therefore provide a different source of directional information.

Motivated by this observation, the first main goal of the present work is to design an efficient and computationally cheap diagonal scaling mechanism for noisy large-scale DFO problems. To this end, we construct diagonal curvature approximations from conjugacy-type relations among generated search directions, rather than from explicit gradient estimates or secant equations involving gradient differences. The resulting update is derivative-free in a strict sense: it requires neither exact gradients nor approximate gradients, and it avoids finite-difference quotients, whose accuracy may deteriorate significantly in the presence of noise.

The second main goal of the present work is to incorporate the proposed diagonal scaling technique into the \texttt{MADFO} framework \cite{MADFO}, with the aim of improving its suitability for large-scale noisy DFO problems. The original \texttt{MADFO} algorithm is a matrix-adaptation evolution strategy that combines mutation, selection, and recombination phases with adaptive mechanisms for generating effective mutation step sizes, together with a randomized nonmonotone line search and heuristic updates of the best point. By embedding the proposed diagonal scaling strategy into this framework, we seek to improve the coordinate-wise adaptation of the mutation process without introducing gradient approximations or expensive dense matrix operations. This is particularly important in large-scale noisy settings, where full matrix adaptation can become computationally demanding, while scalar step-size adaptation alone may fail to capture differences in variable sensitivities. The resulting algorithm preserves the derivative-free and noise-oriented nature of \texttt{MADFO}, while equipping it with a low-cost mechanism for exploiting directional curvature information through diagonal scaling.
The proposed mechanism is also computationally attractive for high-dimensional problems. Since the approximation is diagonal, it requires only vector storage, componentwise operations, and a small number of inner products. Consequently, both storage and arithmetic costs scale linearly with the problem dimension. This makes the approach substantially cheaper than dense matrix-based scaling strategies and suitable for incorporation into large-scale mutation-based DFO frameworks, including \texttt{MAES}-type algorithms and advanced variants such as \texttt{MADFO}. In this way, the proposed method provides a noise-aware and low-cost alternative for constructing diagonal scaling information in high-dimensional DFO.

In the numerical comparisons reported in this paper, robustness and efficiency are interpreted in the standard performance-testing sense. A solver is regarded as more robust if it solves a larger number of test problems within the prescribed computational budget and stopping criteria. Efficiency is measured with respect to the number of objective-function evaluations required to solve a problem; thus, among solvers that solve a given problem, the more efficient solver is the one requiring fewer function evaluations. Accordingly, the reported robustness and efficiency comparisons below refer to these two criteria.

The main contributions of this work are summarized as follows:
\begin{itemize}
\item[--] We incorporate the proposed {\bf diagonal scaling mechanism} into the {\bf mutation phase} of a basic {\tt MAES} framework. The resulting mutation scheme replaces the dense affine scaling matrix by a diagonal inverse-square-root scaling matrix, while preserving the derivative-free and noise-aware structure of the underlying evolutionary framework. This generic construction is then specialized within the {\tt MADFO} framework to obtain the proposed large-scale variant.

\item[--] We derive several derivative-free diagonal curvature-surrogate updates from
conjugacy-type conditions among consecutive generated directions. In particular,
we consider {\bf quadratic-penalty}, {\bf cubic-penalty}, and more
general {\bf higher-order penalty} models, all of which reduce to low-dimensional
or scalar computations.

\item[--] We incorporate the proposed {\bf diagonal scaling mechanism} into the
{\bf mutation phase} of {\tt MAES}. The resulting variant replaces the dense
affine scaling matrix by a diagonal inverse-square-root scaling matrix, while
preserving the derivative-free and noise-aware structure of the underlying
evolutionary framework.

\item[--] We show that the proposed mutation-scaling mechanism requires only
{\bf vector storage}, componentwise operations, and a small number of inner products.
Hence, its storage and arithmetic costs scale {\bf linearly} with the problem
dimension, making it suitable for high-dimensional noisy DFO problems.

\item[--] We implement the proposed mechanism as a large-scale diagonal variant of
{\tt MADFO}, called {\tt DMADFO}, and compare it with {\tt MADFO}
\cite{MADFO} and {\tt LMMAES} \cite{Loshchilov2019} on
noisy DFO test functions obtained from the \texttt{prince} test problems of the {\tt BARON} software \cite{BARON}. On small- to medium-scale noisy DFO problems with dimensions \(2\leq n\leq100\),
the proposed diagonal scaling preserves the robustness and efficiency of full-matrix
{\tt MADFO} up to a loss of less than {\bf \(3\%\)}. On noisy DFO problems ranging
from low to high dimension, with \(2\leq n\leq10000\), {\tt DMADFO} improves
robustness and efficiency by more than {\bf \(20\%\)} compared with {\tt LMMAES}.
These results indicate that the diagonal scaling preserves most of the robustness and
efficiency of full affine scaling while providing a more practical alternative for
high-dimensional noisy DFO.
\end{itemize}



The remainder of this paper is organized as follows. In Section
\ref{diagonal}, we present derivative-free diagonal approximation models
for the Hessian within a least-change penalty framework that not only
promotes conjugacy but also enhances well-conditioning. In Section \ref{MAES}, we show how the proposed low-cost diagonal curvature approximations can be incorporated into the mutation phase of a basic {\tt MAES} framework, and we explain how this idea is specialized within {\tt MADFO} to obtain the proposed {\tt DMADFO} variant. In
Section \ref{Numerical}, we evaluate the proposed approach on noisy benchmark problems, and report the numerical results in detail. Finally,
we summarize the main findings and conclusions of this study in Section
\ref{remarks}.

\section{Diagonal Hessian Approximation via Conjugacy Condition}\label{diagonal}

To the best of our knowledge, the diagonal Hessian approximations proposed in the literature have mainly been developed based on variants of the secant equation or classical finite-difference (derivative) approximations, often within a least-change framework; namely, the new (inverse) Hessian approximation is obtained through the smallest possible modification of the current approximation.  As recent examples, Enshaei et al. \cite{Enshaei2016} employed a least-change model based on a generalized Frobenius norm that incorporates a weak secant equation as a constraint to derive diagonal Hessian approximations. This work was subsequently extended to implicitly incorporate higher-order information of the underlying model \cite{Enshaei2018}. Andrei \cite{Andrei2018} also employed an ordinary least-change model in which the diagonal Hessian approximation is computed by incorporating a low-cost ill-conditioning penalty term while enforcing a weak secant equation. Moreover, he used classical finite-difference estimations based on first-order information to develop such diagonal approximations \cite{Andrei2020,Andrei20202}. In another study within this category, Babaie--Kafaki et al. \cite{BabaieKafaki2022} employed a diagonal least-squares approximation of the secant equation to utilize first-order information from the underlying model, while incorporating the Byrd--Nocedal measure function to promote well-conditioning, as previously done in \cite{Andrei20202}. In addition, Aminifard and Babaie--Kafaki \cite{Aminifard2023} developed diagonal Hessian approximations based on an eigenvalue analysis of the well-known BFGS and DFP quasi--Newton updating formulas, which involve the successive gradient displacement vector as well.

Although such studies have demonstrated promising improvements from a computational perspective, they generally require explicit first-order information  of the relevant NLP model. Therefore, they cannot be directly employed within DFO frameworks. Hence, to make such sparse approximations practical for DFO, it is necessary to incorporate an alternative theoretical feature. It should be emphasized that, in the context of DFO, only objective function values and the search directions (steps) generated by a prescribed rule are available. Accordingly, the geometric structure of the solution trajectory should be carefully analyzed in order to extract useful information for subsequent iterations.

It is widely accepted that the concept of conjugacy is a fundamental notion in linear algebra and numerical optimization \cite{Watkins}. As an extension of orthogonality, it has revolutionized the algorithmic aspects of NLP, mainly through the development of conjugate gradient algorithms as efficient classical approaches for solving large-scale NLP problems and high-dimensional linear systems \cite{SBKRairo}. In particular, it is theoretically established that when the objective function in \eqref{UP} is a strictly convex quadratic function, conjugate gradient algorithms satisfy the quadratic ($n$-step) termination property \cite{SunYuan}. Moreover, their convergence behavior is improved significantly when the eigenvalues of the Hessian matrix associated with \eqref{UP} exhibit low diversity \cite{Watkins}. On the other hand, for such objective functions, it has also been shown that well-known quasi--Newton methods exhibit similar convergence behavior \cite{SunYuan}. Motivated by these observations, the conjugacy condition can be regarded as a natural and constructive underlying property for developing diagonal Hessian approximations within DFO frameworks.

Several fundamental conjugacy conditions have been introduced in the literature \cite{SBKRairo}, where at least first-order information of the underlying model is required to formulate the condition. In order to develop conjugate gradient algorithms for general unconstrained optimization problems, where unlike quadratic models the Hessian matrix is not fixed, conjugacy is typically imposed on successive search directions \cite{DaiLiaoNCG}. Such an approach provides a meaningful framework that can be beneficial for constructing sparse Hessian approximations using only zero-order information of the model. Therefore, in this work, we employ this successive imposition of conjugacy as an alternative to secant equations or finite-difference techniques that require first-order information of the underlying model \cite{Andreibook}.

To describe our estimation approach in detail, assume that
\begin{equation}\label{Pkdiag}
    \mathcal{P}_k
    :=
    \text{diag}(p_{1_k},p_{2_k},\ldots,p_{n_k})
    \approx
    \nabla^2 f(\bx_k)
\end{equation}
is the available diagonal Hessian approximation, with $\mathcal{P}_0=\mathbf{I}$, which should be updated to
\begin{equation*}
    \mathcal{P}_{k+1}:=\text{diag}(p_{1_{k+1}},p_{2_{k+1}},\ldots,p_{n_{k+1}})\approx\nabla^2f(\bx_{k+1}).
\end{equation*}
Now, instead of relying on previously proposed first-order strategies, particularly those based on the classical quasi--Newton framework derived from the secant equation, we require $\mathcal{P}_{k+1}$ to satisfy another fundamental property of quasi--Newton updating formulas, namely the conjugacy condition \cite{SunYuan}. Specifically, when two successive search directions $\bd_{k-1}$ and $\bd_k$ are available, we determine $\mathcal{P}_{k+1}$ such that
\begin{equation}\label{conjugacy}
    \mathcal{C}_{(k-1,k)}(\mathcal{P}_{k+1})
    := \bd_{k-1}^\top\mathcal{P}_{k+1}\bd_k
    = \sum_{i=1}^n \tau_{i}^{(k-1,k)} p_{i_{k+1}}
    =0,
\end{equation}
where 
\begin{equation}\label{tauk}
    \tau_{i}^{(k-1,k)}:=d_{i_{k-1}}d_{i_k}, \qquad i=1,2,\ldots,n,
\end{equation}
in which $d_{i_{k-1}}$ and $d_{i_k}$ are respectively the $i$th entry of the directions $\bd_{k-1}$ and $\bd_k$. 

The conjugacy effect in our approach can be further strengthened by incorporating additional previously generated search directions throughout the underlying algorithm. More precisely, for arbitrary indices $j<l\leq k$, we may also impose
\begin{equation}\label{conjugacy+}
    \mathcal{C}_{(j,l)}(\mathcal{P}_{k+1})
    := \bd_j^\top\mathcal{P}_{k+1}\bd_l
    = \sum_{i=1}^n \tau^{(j,l)}_i p_{i_{k+1}}
    =0,
\end{equation}
where $\tau^{(j,l)}_i:=d_{i_j}d_{i_l}$, $i=1,2,\ldots,n$. By incorporating more available data into the estimation process, the degree of conjugacy can be enhanced. However, the use of excessively old information may deteriorate the accuracy of the approximation. Therefore, within the proposed framework, greater emphasis should be assigned to more recently generated search directions. By the way, from this perspective, the approximation model benefits from a certain degree of flexibility.

As previously mentioned, in many Hessian approximation techniques proposed in the literature, the least-change property plays a central role \cite{Enshaei2016,Andrei2018,Aminifard2023}. By exploiting this property, quasi--Newton updates can also be derived, which is particularly useful for preserving certain information from the previous iteration  \cite{SunYuan}. In the context of our formulation, this property is also incorporated through a least-squares model that minimizes
\begin{equation}\label{phi}
    \phi_k(\mathcal{P}_{k+1})
    := \|\mathcal{P}_{k+1}-\mathcal{P}_k\|_F^2
    = \sum_{i=1}^n (p_{i_{k+1}}-p_{i_k})^2,
\end{equation}
where $\|\cdot\|_F$ denotes the Frobenius norm. Meanwhile, the function is sufficiently flexible to incorporate additional terms that enforce specific theoretical properties, such as well-conditioning \cite{Watkins}, which is also computationally advantageous.

It is worth mentioning that, among the various factors influencing the effectiveness of a Hessian approximation, one of the most important is that the matrix should be properly scaled, or equivalently, well-conditioned \cite{SBKJOTA}. In other words, in the context of our diagonal approximation approach, the distribution of the diagonal entries should not be excessively irregular, in the sense that there should not be a large disparity between the smallest and largest (absolute) diagonal values, which are precisely the eigenvalues of the matrix \cite{SBKRairo}. As is well known, an ill-conditioned matrix approximation may significantly distort the solution trajectory of the underlying algorithm, lead to unreliable approximations, and adversely affect the convergence speed \cite{Watkins}. These observations motivate us to further enhance the conditioning of the resulting diagonal Hessian approximations.

Here, in order to control the conditioning of the successively generated diagonal Hessian approximations, we incorporate an additional penalty term into the approximation models to keep the resulting approximation away from ill-conditioning. At the same time, the solvability of the resulting model must also be taken into account. To this end, we focus on a key factor associated with a well-known family of well-conditioned matrices, namely the class of orthogonal matrices, which may be regarded as optimally conditioned since their spectral condition number attains its minimum possible value, namely one \cite{Watkins}. It is worth noting that the inverse of a symmetric orthogonal matrix (such as a Householder matrix) coincides exactly with the matrix itself \cite{Watkins}. This property has formed the basis of several recent studies \cite{BabaieKafaki2015}, leading to promising computational results as well \cite{SBKRairo}. Motivated by these observations, we also take this feature into account and seek to make $\mathcal{P}_{k+1}$ close not only to $\mathcal{P}_k$, but also to $\mathcal{P}_k^{-1}$, within the framework of the least-change model. Thus, we propose the next conditioning measure:
\begin{equation}\label{psi}
\psi_k (\mathcal{P}_{k+1}):=\|\mathcal{P}_{k+1}-\mathcal{P}_k^{-1}\|_F^2 = \sum_{i=1}^n \left(p_{i_{k+1}}-\dfrac{1}{p_{i_k}}\right)^2,
\end{equation}
which is well-defined whenever the available diagonal Hessian approximation $\mathcal{P}_k$ is nonsingular.


Now, we are in a position to introduce a general estimation framework for constructing diagonal derivative-free Hessian approximations, 
merging the functions defined above. More specifically, the least-change measure $\phi_k(\mathcal{P}_{k+1})$ is adopted as the
principal component of the objective function, while well-conditioning and
conjugacy are enforced through the penalty terms involving
$\psi_k(\mathcal{P}_{k+1})$ and
$\mathcal{C}_{(k-1,k)}(\mathcal{P}_{k+1})$. Accordingly, the proposed estimation model takes the form
\begin{equation}\label{modelgeneral}
\min_{\mathcal{P}_{k+1}\in\mathbb{D}^n}
\quad
\dfrac{1}{2}\phi_k(\mathcal{P}_{k+1})
+
\dfrac{\rho}{2}\psi_k(\mathcal{P}_{k+1})+
\dfrac{\mu}{2p+z}
\left|
\mathcal{C}_{(k-1,k)}(\mathcal{P}_{k+1})
\right|^{2p+z},
\end{equation}
where \(\rho\geq0\) and \(\mu>0\) are penalty parameters, \(p\geq 1\) is an integer, and \(z\in\{0,1\}\). Note that \(z=0\) corresponds to an even-order penalty model and \(z=1\) corresponds to an odd-order penalty model. As can be seen, when $\mu$ is sufficiently large and the Hessian approximation $\mathcal{P}_{k+1}$ fails to satisfy the conjugacy condition adequately, a substantial penalty is imposed on the objective function of \eqref{modelgeneral}. Consequently, the solution process is automatically steered toward approximations that satisfy the conjugacy condition \eqref{conjugacy} more closely, while taking the well-conditioning of the approximations into account through the other penalty term. Meanwhile, in accordance with the results reported in the literature \cite{SunYuan}, the choices of appropriate values for the penalty parameters can significantly influence the quality of the resulting solution. Moreover, larger values of \(p\) yield a flatter response to small conjugacy violations and a steeper response to larger ones. This feature may be particularly beneficial in noisy DFO settings, where minor violations can arise from unreliable directional information, whereas significant violations should still be effectively penalized.

Now, we discuss how the estimation model \eqref{modelgeneral} can be handled for several different choices of $p$ and $z$. The first and foremost observation is that, since $x\mapsto |x|^{2p+z}$ is convex for $2p+z\geq 1$ and
$\mathcal{C}_{(k-1,k)}(\mathcal{P}_{k+1})$ is linear, the model
\eqref{modelgeneral} is convex. Consequently, the first-order
optimality condition is sufficient for global optimality. For notational convenience, let
\begin{equation}\label{pkandtaudef}
\bp_{k+1}:=\begin{pmatrix}p_{1_{k+1}}\\p_{2_{k+1}}\\ \vdots\\p_{n_{k+1}}\end{pmatrix},\quad
 \bp_{k}:=\begin{pmatrix}p_{1_{k}}\\p_{2_{k}}\\ \vdots\\p_{n_{k}}\end{pmatrix},\quad\text{and}\quad
 \btau_k:=\begin{pmatrix}\tau^{(k-1,k)}_1\\\tau^{(k-1,k)}_2\\ \vdots\\\tau^{(k-1,k)}_n\end{pmatrix}\overset{\eqref{tauk}}{=}\bd_{k-1}\odot \bd_k, 
\end{equation}
and also,
\begin{equation}\label{defvecb} 
 \bb_k:=\dfrac{1}{1+\rho}\begin{pmatrix} p_{1_k}+\dfrac{\rho}{p_{1_k}}, p_{2_k}+\dfrac{\rho}{p_{2_k}}, \dots, p_{n_k}+\dfrac{\rho}{p_{n_k}} \end{pmatrix}^\top.
\end{equation}
We also introduce the following diagonal matrices to facilitate a compact and clear presentation of the results:
\begin{eqnarray}
\label{defmatBK}  \mathcal{B}_{k} &:=& \dfrac{1}{1+\rho}\text{diag}\begin{pmatrix} p_{1_k}+\dfrac{\rho}{p_{1_k}}, p_{2_k}+\dfrac{\rho}{p_{2_k}}, \dots, p_{n_k}+\dfrac{\rho}{p_{n_k}} \end{pmatrix},\\
\label{defTk} \mathcal{T}_k &:=& \text{diag}\begin{pmatrix}\tau^{(k-1,k)}_1, \tau^{(k-1,k)}_2,\ldots, \tau^{(k-1,k)}_n\end{pmatrix}.
\end{eqnarray}

\begin{rem}

From a computational perspective, and in order to ensure the well-definedness of the successive approximations, it is necessary to impose appropriate safeguarding conditions on the iterative approximations. In this regard, the first and most important requirement is that all intermediate Hessian approximations remain nonsingular, so that the function $\psi(\cdot)$, and consequently model \eqref{modelgeneral}, are well-defined. Furthermore, it is often desirable to keep the diagonal entries bounded between two positive constants in order to preserve positive definiteness, a key property expected in a neighborhood of the solution, and to avoid numerical difficulties caused by arising excessively small or large (absolute) values \cite{Watkins}. A detailed discussion of such necessary numerical considerations is deferred to experimental parts of this study.

\end{rem}

\subsection{Quadratic Conjugacy Form with $(p,z)=(1,0)$} 
By such a setting, that is practically advisable according to the norm of the literature, the first-order optimality conditions yield
\begin{equation}\label{system1}
 p_{j_{k+1}} + \dfrac{\mu}{1+\rho} \tau^{(k-1,k)}_j \sum_{i=1}^n \tau^{(k-1,k)}_ip_{i_{k+1}} = \dfrac{1}{1+\rho}\left(p_{j_k}+\dfrac{\rho}{p_{j_k}}\right),\qquad j=1,2,...,n.
\end{equation}
Using the notation introduced previously, the system \eqref{system1} can be expressed in the following matrix form:
\begin{equation}\label{system2}
 \left(\mathbf{I} + \dfrac{\mu}{1+\rho} \btau_k\btau_k^\top \right)\bp_{k+1}=\bb_k.
\end{equation}

It can be seen that the coefficient matrix of system \eqref{system2} is a rank-one update of the identity matrix, a class of matrices that has been extensively studied in the literature, particularly in the context of quasi--Newton updating formulas \cite{SunYuan}. Since the penalty parameters $\mu$ and $\rho$ are nonnegative scalars, the coefficient matrix is positive definite. Therefore, by the Sherman--Morrison formula \cite{SunYuan}, the analytical solution of the system is given by
\begin{equation*} 
 \bp^{\quadc}_{k+1} := \left(\mathbf{I} - \dfrac{\mu\ \btau_k\btau_k^\top}{\mu\ \|\btau_k\|^2 +\rho+1}  \right) \bb_k
 =\bb_k -  \dfrac{\mu\ \btau_k^\top\bb_k}{\mu\ \|\btau_k\|^2 +\rho+1}\btau_k.
\end{equation*}
Hence, a diagonal Hessian approximation is obtained without requiring any derivative information from the underlying model, merely by storing a few vectors and computing a small number of inner products. More exactly, independently of derivative information, the following updating formula is available:
\begin{equation}\label{pk1final2}
    \mathcal{P}^{\quadc}_{k+1}:=\mathcal{B}_{k} - \dfrac{\mu\ \btau_k^\top\bb_k}{\mu\ \|\btau_k\|^2 +\rho+1}\ \mathcal{T}_k\approx\nabla^2f(\bx_{k+1}).
\end{equation}

\subsection{Cubic Conjugacy Form with $(p,z)=(1,1)$} 
As a further variant of the penalized estimation model \eqref{modelgeneral}, one may consider the penalty imposed on the conjugacy violation in a cubic form by setting $(p,z)=(1,1)$. The motivation is to penalize large violations of the conjugacy condition more strongly
while still preserving a cheaply computable structure. Now, by defining 
\begin{equation}\label{gamma}
    \gamma_k:=\btau_k^\top\bp_{k+1},
\end{equation}
the first-order optimality conditions of the model in the vector form can be written as 
\begin{equation}\label{cubicFOC}
    \bp_{k+1} +\dfrac{\mu}{1+\rho} \left|\gamma_k\right|\gamma_k \btau_k = \bb_k.
\end{equation}
Taking inner product of both sides of \eqref{cubicFOC} by $\btau_k^\top$, we obtain
\begin{equation*}
    \left(1 +\dfrac{\mu}{1+\rho}\|\btau_k\|^2 \left|\gamma_k\right|\right) \gamma_k = \bb_k^\top\btau_k.
\end{equation*}
So, if $\bb_k^\top\btau_k=0$, then $\bp_{k+1}^{\cubic}:=\bb_k$; otherwise, the following quadratic equation should be solved:
\begin{equation}\label{equationgamma}
    \eta_k |\gamma_k|\gamma_k + \gamma_k - (\bb_k^\top\btau_k) =0,
\end{equation}
in which $\eta_k:=\mu(1+\rho)^{-1}\|\btau_k\|^2$.

Note that the scalar function $\gamma\mapsto \eta_k|\gamma|\gamma+\gamma-\bb_k^\top\btau_k$ is strictly increasing. Therefore, the equation \eqref{equationgamma} has a unique solution, and 
\begin{equation}\label{sign}
    \operatorname{sign}(\gamma_k) = \operatorname{sign}\left(\bb_k^\top\btau_k\right).
\end{equation}
Setting $\Delta_k:= 1 + 4\eta_k|\bb_k^\top\btau_k|$, the solution of \eqref{equationgamma} can be given in the following (numerically stable) form:
\begin{equation}\label{cubicsolution}
    \gamma_k^{\text{sol}} := \text{sign}\left(\bb_k^\top\btau_k\right)\dfrac{2|\bb_k^\top\btau_k|}{1 + \sqrt{\Delta_k}}.
\end{equation}
Hence, the vector of diagonal entries generated by the cubic-penalty model is given by
\begin{equation}\label{pcubicfinal}
    \bp_{k+1}^{\cubic} := \bb_k-\dfrac{\mu}{1+\rho} |\gamma_k^{\text{sol}}|\gamma_k^{\text{sol}}\btau_k,
\end{equation}
where \(\gamma_k^{\text{sol}}\) is given by \eqref{cubicsolution}. Thus, the corresponding diagonal Hessian approximation is defined by
\begin{equation}\label{Pcubicfinal}
    \mathcal{P}_{k+1}^{\cubic}
    :=
    \begin{cases}
    \displaystyle
    \mathcal B_k-\dfrac{\mu}{1+\rho} |\gamma_k^{\text{sol}}|\gamma_k^{\text{sol}} \mathcal T_k,
    &\qquad \text{if } \bb_k^\top\btau_k\neq0,\\[2ex]
    \mathcal B_k,
    &\qquad \text{otherwise}.
    \end{cases}
\end{equation}

\subsection{General Conjugacy Form with $p>1$}

Unlike the previous cases, since now $2p+z>3$, although the penalty associated with conjugacy violation is further strengthened, the corresponding first-order optimality system does not admit an analytical solution. Therefore, one has to resort to iterative numerical procedures to compute an approximate solution. By the way, here the first-order optimality conditions of \eqref{modelgeneral} can be written as
\begin{equation}\label{generalFOC}
    \bp_{k+1}
    +
    \dfrac{\mu}{1+\rho}
    |\gamma_k|^{2p+z-2}\gamma_k\btau_k
    =
    \bb_k,
\end{equation}
with $\gamma_k$ defined by \eqref{gamma}.  Taking the inner product of both sides of \eqref{generalFOC} with \(\btau_k\), we
obtain the scalar equation
\begin{equation}\label{generalScalar}
    \gamma_k
    +
    \dfrac{\mu}{1+\rho}
    \|\btau_k\|^2
    |\gamma_k|^{2p+z-2}\gamma_k
    =
    \bb_k^\top\btau_k,
\end{equation}
in which the left-hand side is strictly
increasing in terms of \(\gamma_k\). Hence, the equation \eqref{generalScalar} has a unique solution and \eqref{sign} holds. So, if \(\bb_k^\top\btau_k=0\), then \(\gamma_k=0\), and therefore $\bp_{k+1}^{(2p+z)}:=\bb_k$; otherwise, the unique solution, denoted by \(\gamma_k^{(2p+z)}\), can be computed by solving a scalar equation. 

Now, considering \eqref{sign}, we write $\gamma_k^{(2p+z)} = \operatorname{sign}\!\left(\bb_k^\top\btau_k\right)u_k$, for some $u_k\geq 0$. Substituting this representation into \eqref{generalScalar} shows that $u_k$ is the unique nonnegative solution of
\begin{equation}\label{LHS}
    u_k
    +
    \dfrac{\mu}{1+\rho}
    \|\btau_k\|^2
    u_k^{2p+z-1}
    =
    \left|\bb_k^\top\btau_k\right|,
\end{equation}

in which, the left-hand side is continuous and strictly increasing for \(u_k\geq0\), with value \(0\) at \(u_k=0\). Moreover, at $u_k=\left|\bb_k^\top\btau_k\right|$, the left-hand side of \eqref{LHS} is at least $\left|\bb_k^\top\btau_k\right|$. Therefore, the solution lies in the interval $\left[0,\left|\bb_k^\top\btau_k\right|\right]$. This observation is sufficient to solve the resulting one-dimensional monotone equation over the bounded interval.

In the implementation, the one-dimensional monotone equation \eqref{LHS} is approximately solved by applying a fixed number of bisection iterations on the interval
\([0,|\bb_k^\top\btau_k|]\). Such a scalar computation is independent of the dimension \(n\), and therefore does not affect the \(O(n)\) cost of the update.  Once \(\gamma_k^{(2p+z)}\) is obtained, the vector of diagonal entries is updated by
\begin{equation}\label{pGeneral}
    \bp_{k+1}^{(2p+z)}
    =
    \bb_k
    -
    \dfrac{\mu}{1+\rho}
    \left|\gamma_k^{(2p+z)}\right|^{2p+z-2}
    \gamma_k^{(2p+z)}
    \btau_k.
\end{equation}
Thus, the corresponding diagonal Hessian approximation can be formulated by
\begin{equation}\label{PGeneral}
\mathcal{P}_{k+1}^{(2p+z)}
=
\begin{cases}
\displaystyle
\mathcal B_k
-
\dfrac{\mu}{1+\rho}
\left|\gamma_k^{(2p+z)}\right|^{2p+z-2}
\gamma_k^{(2p+z)}
\mathcal T_k,
& \qquad\text{if } \bb_k^\top\btau_k \neq 0,\\[2ex]
\mathcal B_k,
& \qquad\text{otherwise},
\end{cases}
\end{equation}
where \(\btau_k\), \(\bb_k\), \(\mathcal B_k\), and \(\mathcal T_k\) are defined in
\eqref{pkandtaudef}, \eqref{defvecb}, \eqref{defmatBK}, and \eqref{defTk},
respectively.


\begin{rem}
    The overall computational cost of the approximation process (in all the above cases) is linear in terms of \(n\), namely \(O(n)\). Indeed, computing \(\btau_k^\top\bb_k\) and \(\|\btau_k\|^2\) requires only \(O(n)\) operations. Solving \eqref{generalScalar} involves a one-dimensional bisection procedure, while forming \(\mathcal{P}_{k+1}^{(2p+z)}\) requires only diagonal matrix operations and therefore also costs \(O(n)\). That is, the proposed approach is not only independent of derivative information, but also provides low-cost, memoryless approximations that can be readily incorporated into DFO algorithms.
\end{rem}

\section{Diagonal Conjugacy-Based Scaling in the Mutation Phase of {\tt MAES}}\label{MAES}

In this section, we describe how the diagonal approximations obtained in the previous section can be incorporated into the mutation phase of {\tt MAES}. The main idea is to replace the affine scaling matrix used to generate mutation directions with a diagonal scaling matrix derived from the proposed derivative-free Hessian approximation. This replacement is particularly motivated by highly noisy regimes, where the sorting and selection phases may identify misleading offspring and, consequently, the recombination point may provide only imperfect directional information. Instead of using this uncertain information to update a dense affine or covariance matrix, the proposed strategy converts it into a conservative diagonal scaling update. Thus, the method does not assume that the generated directions are noise-free; rather, it uses them cautiously, through a bounded and regularized coordinatewise mechanism. In this way, the influence of noisy rankings is restricted, while the derivative-free structure and low computational cost of the evolutionary framework are preserved.

The evaluation strategy is a randomized DFO framework based on repeated interactions among three main phases: mutation, selection, and recombination; see, e.g., \cite{Araki2014}. In the mutation phase, new candidate solutions are generated by perturbing the current parental information, usually through random vectors with zero mean. In the selection phase, the generated candidate solutions are ranked according to their possibly inexact function values, and a subset of the most promising individuals is selected to contribute to the next generation. In the recombination phase, the selected individuals are combined to define the new parental information for the following iteration. This mechanism originates from the principle of biological evolution: mutation generates new individuals, selection identifies promising candidates, and recombination transfers information to the next generation.

A well-known class of randomized DFO methods based on this principle is \texttt{CMAES}, namely the covariance matrix adaptation evolution strategy. These methods are designed for solving nonlinear, nonconvex, and possibly noisy continuous optimization problems. In \texttt{CMAES}, the candidate solutions are sampled from a multivariate normal distribution, whose covariance matrix determines the scaling of the search distribution and the dependencies between the variables. The adaptation of this covariance structure enables the method to learn useful geometric information about the objective function and is especially beneficial for ill-conditioned problems.

Several variants of \texttt{CMAES} have been proposed in the literature, mainly differing in the way the covariance or affine scaling information is updated. Examples include the classical \texttt{CMAES} solver of Auger and Hansen \cite{AugH}, the limited-memory \texttt{LMMAES} solver of Loshchilov et al. \cite{Loshchilov2019}, the fast matrix adaptation strategy \texttt{fMAES} of Beyer \cite{Beyer2020}, the \texttt{BiPopMAES} variant of Beyer and Sendhoff \cite{Beyer2017}, and the {\tt MADFO} and {\tt MATRS} solvers of Kimiaei and Neumaier \cite{MADFO,MATRS}. From a second-order viewpoint, the covariance adaptation mechanism can be interpreted as a derivative-free analogue of constructing a local curvature model. This is conceptually related to quasi--Newton methods, where Hessian or inverse Hessian approximations are used to scale the search directions.

A key advantage of \texttt{CMAES}-type methods is that no explicit derivative information is required. In many implementations, only the rank ordering of the candidate solutions, induced by their possibly inexact function values, is used to update the sampling distribution. Thus, compared with traditional gradient-based optimization schemes, these methods impose fewer smoothness requirements on the objective function and can be applied effectively in noisy derivative-free settings.

Nevertheless, this flexibility also introduces some difficulties. Since the update mechanism is mainly driven by ranking information, the algorithm may update its step-size or scaling information without explicitly verifying whether a reliable decrease in the inexact function value has occurred. In noisy problems, the ordering of the candidate solutions may be distorted by perturbations in the objective values. Consequently, some accepted or highly ranked candidate solutions may still be far from an approximate stationary point of the underlying unconstrained noisy DFO problem, a point whose true gradient norm is below a given threshold. Moreover, maintaining and updating a dense covariance or affine scaling matrix can be computationally demanding in high-dimensional settings.

These observations motivate the use of a diagonal derivative-free scaling mechanism in the mutation phase. The diagonal approximations developed in the previous section provide a low-cost alternative to dense affine scaling. Instead of adapting a full covariance matrix, we construct a diagonal curvature surrogate from the generated search directions and use it to scale the mutation vectors. In this way, the derivative-free character of the evaluation strategy is preserved, while curvature-related information is incorporated through a computationally inexpensive diagonal matrix.

Let \(\bx_k^{\rec}\in\mathbb R^n\) denote the current recombination point at iteration \(k\) (\(\bx_0^{\rec}:=\bx_0\)), and let \(\sigma_k>0\) be the global step-size parameter. In a standard {\tt MAES} framework, the mutation (sample) points are commonly generated as
\begin{equation}\label{standardMAESmutation}
    \bx_k^{(r)}
    =
    \bx_k^{\rec}
    +
    \sigma_k A_k \bz_k^{(r)},
    \qquad
    r=1,2,\ldots,\lambda,
\end{equation}
where \(A_k\in\mathbb R^{n\times n}\) is the affine scaling matrix, \(\lambda\) is the number of mutation points or sample size, and the distribution directions
\[
    \bz_k^{(r)}\sim\mathcal N(\mathbf 0,\mathbf I),
    \qquad
    r=1,2,\ldots,\lambda.
\]
The matrix \(A_k\) controls the shape of the mutation distribution. There are several formulas to approximate \(A_k\), e.g., see \cite{AugH,Beyer2017,Beyer2020}. However, in high-dimensional problems, storing and updating a dense affine matrix can be computationally expensive. Therefore, we propose to replace \(A_k\) by a diagonal scaling matrix derived from the diagonal approximation \(\mathcal P_{k+1}\).

At iteration \(k\), the safeguarded diagonal approximation \(\mathcal P_k\), obtained from the previous iteration, plays the role of a diagonal Hessian approximation or derivative-free curvature surrogate. Therefore, the current mutation covariance is naturally chosen proportional to \(\mathcal P_k^{-1}\). Hence, we define
\begin{equation}\label{DkdefMAES}
    D_{k}:=\mathcal P_{k}^{-1/2}=\operatorname{diag}
    \left(
        \frac{1}{\sqrt{p_{1_k}}},
        \frac{1}{\sqrt{p_{2_k}}},
        \ldots,
        \frac{1}{\sqrt{p_{n_k}}}
    \right).
\end{equation}
Therefore, the mutation direction is generated by
\begin{equation}\label{mutationdirectionMAES}
    \by_k^{(r)}
    =
    D_{k}\bz_k^{(r)},
    \qquad
    \bz_k^{(r)}\sim\mathcal N(\mathbf 0,\mathbf I),
\end{equation}
and the mutation points are generated as
\begin{equation}\label{offspringMAESdiagonal}
    \bx_k^{(r)}
    =
    \bx_k^{\rec}
    +
    \sigma_k \by_k^{(r)},
    \qquad
    r=1,2,\ldots,\lambda.
\end{equation}
Thus,
\[
    \by_k^{(r)}
    \sim
    \mathcal N(\mathbf 0,D_{k}D_{k}^{\top})
    =
    \mathcal N(\mathbf 0,\mathcal P_{k}^{-1}).
\]
This means that mutation is enlarged along coordinates with smaller estimated curvature and reduced along coordinates with larger estimated curvature.

To construct \(\mathcal P_{k+1}\), one may use the quadratic penalized conjugacy-based approximation \eqref{pk1final2}. This approximation usually gives a smoother scaling update. Another possibility is to use the cubic-penalty conjugacy model \eqref{Pcubicfinal}. This variant penalizes large conjugacy violations more strongly while producing milder corrections when the scalar conjugacy violation is small. Therefore, the cubic model may be useful in noisy settings, where small violations may be caused by unreliable direction information, whereas large violations should still be discouraged.

It remains to define the search directions used in \eqref{pkandtaudef}, i.e., $\btau_k=\bd_{k-1}\odot\bd_k$. Within {\tt MAES}, a natural derivative-free choice is the normalized displacement of the recombination point. After evaluating and ranking the mutation points, the recombination point is updated by
\begin{equation}\label{meanupdateMAES}
    \bx_{k+1}^{\rec}
    =
    \sum_{r=1}^{\mu_{\rm sel}}
    w_r
    \bx_k^{(r:\lambda)},
\end{equation}
where \(\mu_{\rm sel}\) is the number of selected offspring, \(w_r>0\) are recombination weights satisfying
\[
    \sum_{r=1}^{\mu_{\rm sel}}w_r=1,
\]
and \(\bx_k^{(r:\lambda)}\) denotes the \(r\)th best mutation points among the \(\lambda\) generated candidates. Then, we define the generated displacement direction
\begin{equation}\label{dkMAES}
    \bd_k
    :=
    \frac{\bx_{k+1}^{\rec}-\bx_k^{\rec}}{\sigma_k}.
\end{equation}
This direction reflects the aggregate recombination displacement produced by mutation, selection, and recombination.

Since the next scaling matrix \(D_{k+1}=\mathcal P_{k+1}^{-1/2}\) will contain inverse square roots of the diagonal entries of \(\mathcal P_{k+1}\), positivity must be enforced when constructing \(\mathcal P_{k+1}\). Therefore, after computing a raw diagonal vector \(\widehat{\bp}_{k+1}\), we apply the safeguarding rule
\begin{equation}\label{safeguardMAES}
    p_{i_{k+1}}
    =
    \min
    \left\{
        p_{\max},
        \max
        \left\{
            p_{\min},
            \widehat p_{i_{k+1}}
        \right\}
    \right\},
    \qquad
    i=1,2,\ldots,n,
\end{equation}
where
\[
    0<p_{\min}<p_{\max}<\infty.
\]
This guarantees that \(\mathcal P_{k+1}\) remains nonsingular and positive definite, and that \(D_{k+1}\) is well-defined.

\paragraph{Step-size update.}
The scalar parameter \(\sigma_k>0\) is kept as a global mutation step-size, while the diagonal matrix \(D_k\) controls the coordinatewise scaling of the current mutation distribution. Thus, the roles of \(\sigma_k\) and \(D_k\) are separated: \(D_k\) determines the shape of the current diagonal search distribution, whereas \(\sigma_k\) determines its overall radius. This is consistent with the usual parameterization of evolution strategies, where the global step-size and the covariance or transformation matrix are adapted separately. In the original \texttt{MADFO} solver, the step-size is updated by a traditional \texttt{MAES}-type formula, but it is also adjusted by the function \texttt{goodStepSize} so that the mutation and recombination step sizes are neither too small nor too large \cite{MADFO}. This safeguard is important in noisy problems, because excessively small step sizes may lead to null steps, whereas excessively large step sizes may lead to line-search failures or poorly distributed mutation points.

In the present diagonal variant, the same principle can be retained. After the offspring have been generated according to \eqref{offspringMAESdiagonal}, the recombination point is computed from the selected candidates, and the generated displacement direction is normalized by the step-size that was actually used to produce the mutation and recombination points. Hence, the normalization in the definition of \(\bd_{k}\) uses \(\sigma_k\), not \(\sigma_{k+1}\). The next global step-size \(\sigma_{k+1}\) is updated only after the recombination phase. For example, using a cumulative step-size adaptation mechanism, one may define the recombination distribution direction by
\begin{equation}\label{zrecombMAES}
    \bz_k^{\rm rec}
    :=
    \sum_{r=1}^{\mu_{\rm sel}}
    w_r \bz_k^{(r:\lambda)},
\end{equation}
where \(\mu_{\rm sel}\) is the number of selected candidates, \(w_r>0\) are recombination weights satisfying
\[
    \sum_{r=1}^{\mu_{\rm sel}}w_r=1,
\]
and \(\bz_k^{(r:\lambda)}\) denotes the distribution direction corresponding to the \(r\)th best mutation points. The evolution path for the step-size is then updated by
\begin{equation}\label{psigmaupdateMAES}
    \bp_{\sigma,k+1}
    =
    (1-c_\sigma)\bp_{\sigma,k}
    +
    c_\sigma \bz_k^{\rm rec},
\end{equation}
where \(\bp_{\sigma,k}\in\mathbb R^n\) is the cumulative step-size path and \(c_\sigma\in(0,1]\) is the cumulation parameter. As in \cite{Beyer2020}, the global step-size may then be updated by
\begin{equation}\label{sigmaupdate}
    \sigma_{k+1}
    =
    \sigma_k
    \exp
    \left(
        \frac{c_\sigma}{d_\sigma}
        \left(
            \frac{\|\bp_{\sigma,k+1}\|}{e_\sigma}
            -1
        \right)
    \right),
\end{equation}
where \(d_\sigma>0\) is a damping parameter and \(e_\sigma>0\) is an approximation of the expected norm of a standard Gaussian vector in \(\mathbb R^n\); that is,
\[
    e_\sigma \approx \mathbb E\| \bu \|,
    \qquad
    \bu\sim\mathcal N(\mathbf 0,\mathbf I).
\]
Finally, to avoid excessively small or excessively large global steps, the updated step-size can be projected onto a prescribed interval
\[
    0<\sigma_{\min}<\sigma_{\max}<\infty,
\]
by setting
\begin{equation}\label{sigmasafeguard}
    \sigma_{k+1}
    :=
    \min\{\sigma_{\max},\max\{\sigma_{\min},\sigma_{k+1}\}\}.
\end{equation}
Since the diagonal scaling matrix already performs coordinatewise scaling, this step-size update should be interpreted only as a global radius-control mechanism, not as a replacement for the diagonal curvature scaling.

We define the candidate raw diagonal update by choosing one of the three available
diagonal approximation formulas:
\begin{equation}\label{rawpchoice}
{\mathcal P}^{\text{alg}}_{k+1}
\in
\left\{
\mathcal P^{\quadc}_{k+1},
\mathcal P_{k+1}^{\cubic}, \mathcal{P}_{k+1}^{(2p+z)}
\right\},
\end{equation}
where \(\mathcal P^{\quadc}_{k+1}\) is obtained from the quadratic penalized model
\eqref{pk1final2},  \(\mathcal P_{k+1}^{\cubic}\) is obtained from the cubic-penalty
model \eqref{Pcubicfinal}, and $\mathcal{P}_{k+1}^{(2p+z)}$ by \eqref{PGeneral} with integer \(p\ge 2\) and \(z\in\{0,1\}\).

The resulting diagonal conjugacy-based mutation phase can be summarized as Algorithm \ref{a.DMADFO}, below. A detailed comparison of the storage and arithmetic costs of the proposed diagonal mutation scaling with existing matrix-adaptation strategies is given in Section 2 of the supplemental material \cite{suppMat} of the present paper.

\begin{algorithm}[h]
\caption{The \(k\)th iteration of {\tt MAES} with diagonal conjugacy-based scaling}\label{a.DMADFO}
\begin{algorithmic}[1]
\Require Mean vector \(\bx_k^{\rec}\), step-size \(\sigma_k\), safeguarded diagonal approximation \(\mathcal P_k=\operatorname{diag}(\bp_k)\) obtained from the previous iteration, previous direction \(\bd_{k-1}\), parameters \(\mu>0\), \(\rho\ge0\), bounds \(0<p_{\min}<p_{\max}\), number of mutation points \(\lambda\), number of selected points \(\mu_{\rm sel}\), and recombination weights \(w_r\)

\State Compute the diagonal scaling matrix \(D_k:=\mathcal P_k^{-1/2}\).

\For{\(r=1,2,\ldots,\lambda\)}
    \State Generate the distribution direction \(\bz_k^{(r)}\sim\mathcal N(\mathbf 0,\mathbf I)\).
    \State Compute the mutation direction \(\by_k^{(r)}:=D_k\bz_k^{(r)}\).
    \State Compute the mutation point \(\bx_k^{(r)}\) by \eqref{offspringMAESdiagonal}.
\EndFor

\State Sort the mutation points according to their objective-function values and select the \(\mu_{\rm sel}\) best mutation points
\[
\left\{\bx_k^{(r:\lambda)}\right\}_{r=1}^{\mu_{\rm sel}}.
\]

\State Compute the recombination point \(\bx_{k+1}^{\rec}\) by \eqref{meanupdateMAES}.

\State Compute the generated direction \(\bd_k\) by \eqref{dkMAES}.

\State Compute \(\btau_k\) by \eqref{pkandtaudef} and \(\mathcal B_k\) by \eqref{defmatBK}. Then, choose \({\mathcal P}^{\text{alg}}_{k+1}\) from \eqref{rawpchoice}.

\State Apply the componentwise safeguard \eqref{safeguardMAES} to \({\mathcal P}^{\text{alg}}_{k+1}\) and set \(\mathcal P_{k+1}:={\mathcal P}^{\text{alg}}_{k+1}\).

\State Compute the recombination distribution direction \(\bz_k^{\rm rec}\) by \eqref{zrecombMAES}.

\State Update the evolution path \(\bp_{\sigma,k+1}\) by \eqref{psigmaupdateMAES}.

\State Update the step-size \(\sigma_{k+1}\) by \eqref{sigmaupdate}, and then apply the safeguard \eqref{sigmasafeguard}.
\end{algorithmic} 
\end{algorithm}

\section{Numerical Experiments}\label{Numerical}

In this section, we assess the numerical behavior of the proposed large-scale diagonal variant of \texttt{MADFO}, denoted by \texttt{DMADFO}. The method is based on replacing the full affine scaling matrix used in \texttt{MADFO} by a diagonal inverse-square-root scaling matrix, while preserving the remaining components of \texttt{MADFO}, including selection, recombination, the randomized nonmonotone line search, and the step-size adaptation mechanism. The diagonal scaling is constructed from derivative-free conjugacy-based diagonal Hessian approximations, which exploit only successive search directions and do not require gradient or subgradient information. We also compare a selected variant of \texttt{DMADFO} with \texttt{LMMAES}, which employs a limited-memory covariance matrix approximation, in order to assess whether the proposed diagonal approximation remains competitive against a memory-efficient covariance-based method on high-dimensional noisy problems. This comparison protocol is also motivated by the extensive numerical assessment of \texttt{MADFO} reported in \cite{MADFO}, where \texttt{MADFO} was compared with a large collection of state-of-the-art noisy DFO solvers on a broad set of test problems. Therefore, in the present study, we focus on comparing \texttt{DMADFO} with its full-matrix counterpart \texttt{MADFO}, in order to isolate the effect of the proposed diagonal scaling, and with \texttt{LMMAES}, as a representative limited-memory solver for high-dimensional noisy optimization.

We consider twelve diagonal variants of \texttt{DMADFO}, all based on the general
odd- and even-order penalty model in \eqref{PGeneral}. More precisely,
\(\texttt{DMADFO1}, \ldots, \texttt{DMADFO10}\) correspond to
\(p=1,\ldots,10\), respectively, while \texttt{DMADFO11} and \texttt{DMADFO12}
correspond to \(p=15\) and \(p=20\), respectively. In each case, the affine
scaling matrix used in the mutation phase of \texttt{MADFO} is replaced by the
diagonal scaling matrix generated from \eqref{PGeneral}. Thus, the resulting
variants differ only in the order of the penalty imposed on the conjugacy
violation. Since the update requires only vector operations, a small number of
inner products, and a scalar one-dimensional solve, all variants preserve linear
storage and arithmetic cost with respect to the problem dimension.

A preliminary comparison among the tested \texttt{DMADFO} variants is reported in
Section~3 of the supplementary material \cite{suppMat}. Based on this comparison,
\texttt{DMADFO12} is selected as the representative implementation of
\texttt{DMADFO} for the subsequent comparisons with competing solvers.

We first compare \texttt{DMADFO12} with the original \texttt{MADFO} solver on
small- to medium-scale noisy DFO problems with dimensions \(2\leq n\leq100\).
This comparison is used to assess whether the proposed diagonal scaling can
preserve the robustness of the full affine scaling matrix used in \texttt{MADFO}.
We then compare \texttt{DMADFO12} with \texttt{LMMAES} on low- to large-scale
noisy DFO problems with dimensions \(1\leq n\leq10000\). This second comparison
evaluates whether the selected diagonal variant remains competitive with a
limited-memory affine scaling strategy while being computationally cheaper and
more suitable for high-dimensional noisy settings.

{\bf Stopping Criteria.} The convergence behavior of a solver \(s\) toward a minimizer of the underlying smooth true function \(f\) is assessed through the ratio
\lbeq{e.qs}
q_s:=\frac{f_s-f_{\opt}}{f_0-f_{\opt}},
\qquad s\in\mathcal S,
\eeq
where \(\mathcal S\) denotes the set of solvers under comparison. Here, \(f_s\) is the best objective value obtained by solver \(s\), \(f_0\) is the objective value at the common starting point, and \(f_{\opt}\) is the value at the best known point. The latter is usually associated with a global minimizer, or at least with a high-quality local minimizer, obtained by applying a collection of local and global derivative-free solvers. It should be emphasized that the quantity \(q_s\) is used only for benchmarking purposes, since \(f_{\opt}\) is generally unavailable in real-life applications.

A problem is declared \textit{solved} by solver \(s\) if
\[
    q_s\leq \eps,
\]
and neither the maximum number of function evaluations, denoted by \texttt{nfmax}, nor the maximum allowed CPU time, denoted by \texttt{secmax}, is reached. Otherwise, the problem is declared \textit{unsolved}. The parameters \(\eps\), \texttt{secmax}, and \texttt{nfmax} are chosen so that the best-performing solver is able to solve at least half of the test problems, except in cases where the noise level is so high that increasing \texttt{secmax} or \texttt{nfmax} no longer improves the observed efficiency or robustness. In our experiments, the following choices were found to be effective:
\[
    \texttt{secmax}=600,
    \qquad
    \texttt{nfmax}=10000,
    \qquad
    \eps=10^{-4},
    \qquad
    1\leq n\leq10000.
\]

\bfi{Noisy Test Functions.} We consider four noise models, namely absolute uniform, absolute Gaussian,
relative uniform, and relative Gaussian noises. For each noise model, six
noise levels are tested:
\[
\omega \in \left\{10^{-3},10^{-2},10^{-1},1,10,100\right\}.
\]

For the comparison between the different \texttt{DMADFO} variants and
\texttt{MADFO}, we use the 655 test functions from the \texttt{prince} test problems of the {\tt BARON} software \cite{BARON}
collection with dimensions satisfying \(2\leq n\leq100\), which is available at

``\url{https://minlp.com/optimization-test-problems}".

Each problem is
extended to the four noise models and six noise levels described above.
Consequently, the complete test set consists of
\[
4 \times 6 \times 655 = 15720
\]
noisy test instances.

For the comparison between the most robust \texttt{DMADFO} variant and
\texttt{LMMAES}, we use 926 problems from the \texttt{prince}
collection with dimensions ranging from \(1\) to \(10000\). As in the
previous experiments, each problem is considered under the same four
noise models and six noise levels. Therefore, the complete test set
contains
\[
4 \times 6 \times 926 = 22224
\]
noisy test instances.

\bfi{Initial Points.} Following \cite{KimN,MADFO}, the initial point \(y^0:=0\) is shifted by the vector \(\xi\in\mathbb R^n\), whose components are defined by
\[
    \xi_i
    :=
    (-1)^{i-1}\frac{2}{2+i},
    \qquad
    i=1,\ldots,n.
\]
This shift is applied because some toy problems in the {\tt prince} collection
may have simple solutions that a solver could easily infer from the original
starting point. By shifting the starting point, we reduce the risk of such trivial guessability. Therefore, we choose the shifted initial point
\[
    y^0:=\xi,
\]
and define the initial inexact function value by
\[
    \wt f_0:=\wt f(y^0).
\]
For subsequent iterations, the inexact function values are computed according to
\[
    \wt f_\ell:=\wt f(y^\ell+\xi),
    \qquad
    \ell\geq0.
\]

{\bf Data and Performance Profiles.} To assess the \textit{robustness} and \textit{efficiency} of the compared solvers, we use the data profiles introduced by Mor\'e and Wild \cite{MorW} and the performance profiles proposed by Dolan and Mor\'e \cite{DolM}. Let \(\mathcal S\) be the set of solvers under comparison, and let \(\mathcal P\) be the set of test problems. For a solver \(s\in\mathcal S\), the data profile measures the proportion of problems that can be solved within a budget of \(\kappa\) groups of \(n_p+1\) function evaluations. It is defined as
\lbeq{e.datap}
\delta_{s}(\kappa):=\frac{1}{|\mathcal P|}\Big|\Big\{p\in\mathcal P~ \Big|~cr_{p,s}:=\D\frac{c_{p,s}}{n_{p}+1}\le \kappa\Big\}\Big|.
\eeq
Here, \(n_p\) denotes the dimension of problem \(p\), \(c_{p,s}\) is the \textit{cost measure} required by solver \(s\) to solve problem \(p\), and \(cr_{p,s}\) denotes the corresponding \textit{cost ratio}.

The performance profile of solver \(s\) is given by
\lbeq{e.perp}
\rho_{s}(\tau):=\frac{1}{|\mathcal P|}\Big|\Big\{p\in\mathcal P ~\Big|~ pr_{p,s}:=\D\frac{c_{p,s}}{\min(c_{p,\ol {s}}\mid \ol{s}\in \mathcal S)}\le \tau\Big\}\Big|.
\eeq
This profile represents the fraction of problems for which the performance ratio \(pr_{p,s}\) does not exceed \(\tau\). In particular, \(\rho_s(1)\) gives the fraction of problems on which solver \(s\) performs best relative to the other solvers. For sufficiently large values of \(\tau\) or \(\kappa\), the values of \(\rho_s(\tau)\) and \(\delta_s(\kappa)\), respectively, indicate the proportion of problems solved by solver \(s\).


\subsection{Results on Noisy Smooth DFO problems}

Here, we describe our experimental results on noisy problems, divided into the following two categories, where the proposed algorithms are compared separately with \texttt{MADFO} and \texttt{LMMAES}.

\subsubsection{Comparison of the Most Robust \texttt{DMADFO} Variant with \texttt{MADFO}}

Having identified \texttt{DMADFO12} as the most robust variant of
\texttt{DMADFO}, we now compare it with the original \texttt{MADFO}
algorithm on the same set of $15720$ noisy test functions. Recall that
these instances are obtained from the 655 problems of the
\texttt{prince} collection by considering four noise models
(absolute/relative uniform and Gaussian noises) and six noise levels
$\omega \in \{10^{-3},10^{-2},10^{-1},1,10,100\}$. The purpose of this experiment is to demonstrate that \texttt{DMADFO}, despite
using only a diagonal approximate covariance matrix, remains competitive with
\texttt{MADFO}, which employs a full approximate covariance matrix, across a
broad range of noisy optimization scenarios.

Tables~4 and~5 in \cite{suppMat}
compare the performance of \texttt{MADFO} and the most robust variant of
\texttt{DMADFO}, namely \texttt{DMADFO12}, on the same set of
$15720$ noisy test instances. Overall, \texttt{MADFO} successfully solves
$11112$ problems, whereas \texttt{DMADFO12} solves $10835$ problems.
Thus, \texttt{DMADFO12} solves only 2.49\% fewer problems than
\texttt{MADFO} and retains 97.51\% of its solving capability. Furthermore, the difference in the numbers of unsolved problems is
entirely due to the \texttt{nfmax} termination criterion, while neither
solver experiences failures caused by \texttt{secmax} limits or
algorithmic breakdowns.

These results are particularly encouraging because \texttt{MADFO}
employs a full approximate covariance matrix, whereas
\texttt{DMADFO12} uses only a diagonal approximation. Since a full
covariance matrix can capture correlations between variables, one
might expect a noticeable deterioration in performance when replacing it
with a diagonal approximation. However, the numerical results indicate
that the reduction in robustness is only 2.49\%, despite the
substantial simplification of the covariance model. Therefore, the
diagonal approximation used in \texttt{DMADFO12} appears to capture most
of the relevant scaling information while significantly reducing the
computational and storage requirements associated with maintaining a
full covariance matrix. These results demonstrate that
\texttt{DMADFO12} achieves robustness comparable to that of
\texttt{MADFO}, thereby validating the use of a diagonal approximate
covariance matrix as an efficient alternative to a full approximate
covariance matrix.

\begin{figure}[H]
\centering
\scalebox{0.5}{\includegraphics[width=15cm,height=10cm]{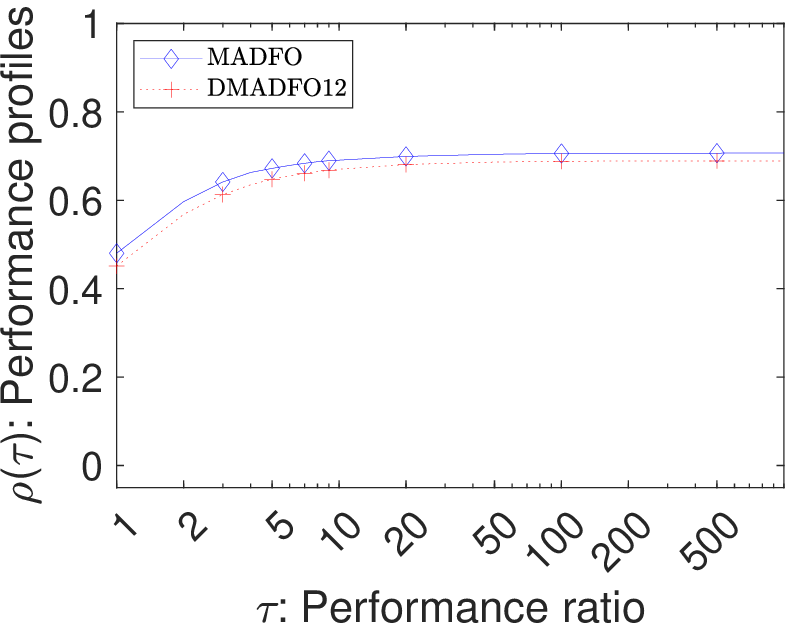}\hspace{0.5cm}\includegraphics[width=15cm,height=10cm]{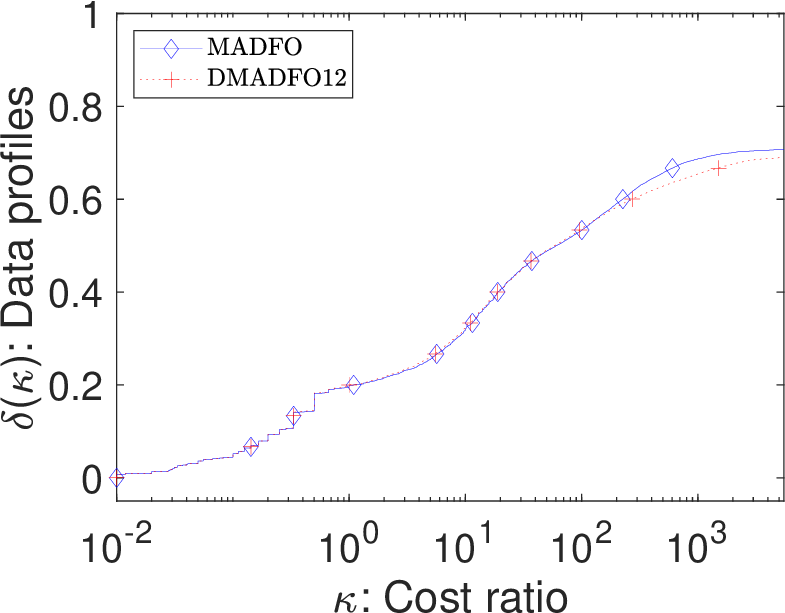}}
\caption{Performance and data profiles of \texttt{MADFO} and \texttt{DMADFO12} based on the number of function evaluations (\texttt{nf}) over all noisy test instances.}
\label{f.f1}
\end{figure}

Figure~\ref{f.f1} presents the performance and data profiles of
\texttt{MADFO} and \texttt{DMADFO12} with respect to the number of
function evaluations. The performance profiles show that
\texttt{DMADFO12}, despite using only a diagonal approximate covariance
matrix, achieves performance close to that of \texttt{MADFO}, with an
efficiency reduction of less than $5\%$ in terms of function
evaluations. This observation further supports the conclusion that most
of the useful scaling information is captured by the diagonal
approximation while avoiding the computational and storage costs
associated with maintaining a full approximate covariance matrix.
Consistent with the results reported in
Tables~4 and 5 in \cite{suppMat}, the
data profiles in Figure~\ref{f.f1} indicate only a minor reduction in
robustness for \texttt{DMADFO12} relative to \texttt{MADFO}.

\subsubsection{Comparison of the Most Robust \texttt{DMADFO} Variant with \texttt{LMMAES}}

Having identified \texttt{DMADFO12} as the most robust variant of
\texttt{DMADFO}, we now compare it with \texttt{LMMAES}, a DFO method based on a limited-memory approximation of the covariance
matrix. The experiments are conducted on 926 test problems from the
\texttt{prince} collection with dimensions ranging from $2$ to $10000$.

As in the previous experiments, each problem is considered under four noise
models, namely absolute uniform, absolute Gaussian, relative uniform, and
relative Gaussian noises, and six noise levels
\[
\omega \in \left\{10^{-3},10^{-2},10^{-1},1,10,100\right\}.
\]
Consequently, the complete test set consists of
\[
4 \times 6 \times 926 = 22224
\]
noisy test instances.

The purpose of this experiment is to compare the robustness of
\texttt{DMADFO12}, which employs a diagonal approximate covariance matrix,
with that of \texttt{LMMAES}, which uses a limited-memory covariance matrix
approximation. In particular, we are interested in assessing whether the
simple diagonal approximation adopted by \texttt{DMADFO12} remains competitive
on large-scale noisy optimization problems, where the storage and computational
costs associated with richer covariance representations become increasingly
significant.

Tables~6 and~7 in \cite{suppMat}
show that \texttt{DMADFO12} clearly outperforms \texttt{LMMAES} on the
$22224$ noisy test functions. In total, \texttt{DMADFO12} solves
$14760$ problems, whereas \texttt{LMMAES} solves only $7184$ problems.
Thus, \texttt{DMADFO12} solves $7576$ more problems than \texttt{LMMAES},
corresponding to about $105.5\%$ more solved problems. Equivalently,
\texttt{DMADFO12} solves $66.4\%$ of all noisy test instances, while
\texttt{LMMAES} solves only $32.3\%$. The unsolved-problem statistics in
Table~7 in \cite{suppMat} confirm the same conclusion:
\texttt{DMADFO12} terminates due to \texttt{nfmax} in $7464$ cases,
whereas \texttt{LMMAES} reaches \texttt{nfmax} in $15040$ cases. Hence,
the diagonal approximate covariance matrix used by \texttt{DMADFO12}
appears substantially more robust than the limited-memory covariance
approximation used by \texttt{LMMAES} on this large-scale noisy test set.

\begin{figure}[H]
\centering
\scalebox{0.5}{%
\includegraphics[width=15cm,height=10cm]{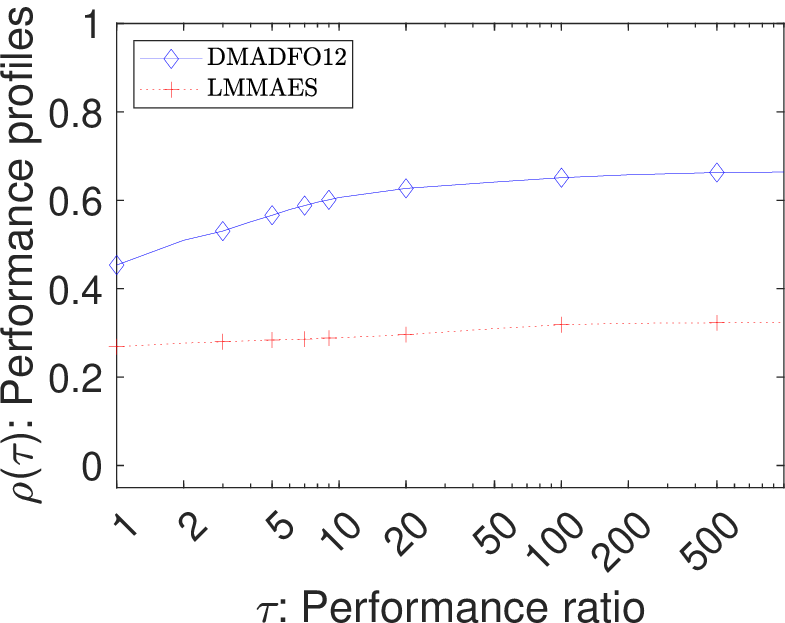}
\hspace{0.5cm}
\includegraphics[width=15cm,height=10cm]{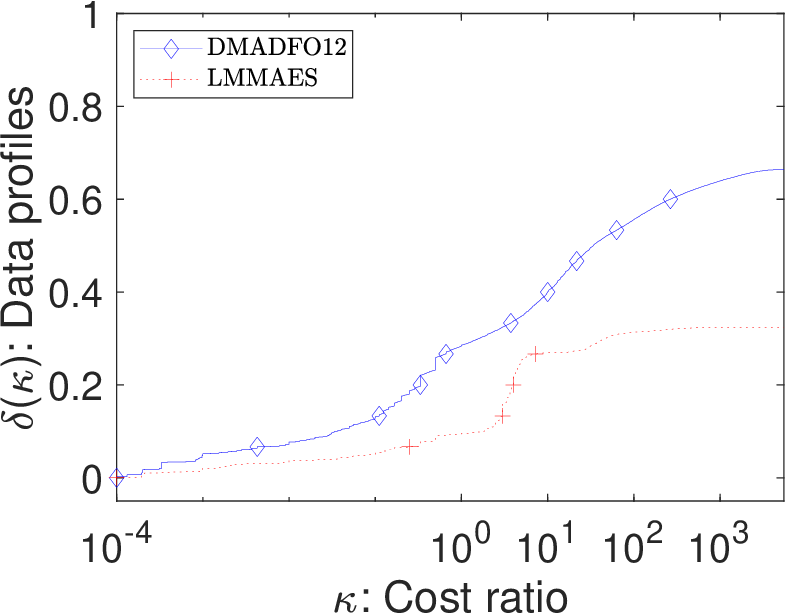}}
\caption{Performance and data profiles of \texttt{LMMAES} and \texttt{DMADFO12} based on the number of function evaluations (\texttt{nf}) over all noisy test instances.}
\label{f.f2}
\end{figure}

Figure~\ref{f.f2} presents the performance and data profiles of
\texttt{LMMAES} and \texttt{DMADFO12} with respect to the number of
function evaluations. The performance profiles show that
\texttt{DMADFO12} is approximately $20\%$ more efficient than
\texttt{LMMAES} in terms of \texttt{nf}-efficiency. This indicates that
the diagonal approximate covariance matrix used in \texttt{DMADFO12}
provides a more effective search behavior than the limited-memory
covariance approximation used in \texttt{LMMAES} on the considered noisy
large-scale test instances. Consistent with
Tables~6 and~7 in \cite{suppMat},
the data profiles in Figure~\ref{f.f2} further confirm the superior
robustness of \texttt{DMADFO12}.

\section{Concluding Remarks}\label{remarks}

We proposed a diagonal affine-scaling derivative-free optimization method for
large-scale noisy problems. The motivation comes from a structural weakness of matrix-adaptation evolution strategies in high-noise
regimes. Although {\tt MAES}-type methods are often more stable than many
traditional derivative-free methods under moderate noise, their performance may
deteriorate when noisy function values corrupt the sorting and selection phases.
In that case, recombination may be based on misleadingly ranked sampled points,
and the information used to update the affine-scaling or matrix-adaptation
mechanism may no longer reflect the local geometry of the objective. This can
make the resulting scaling information inaccurate and reduce the efficiency of
the search.

The proposed method addresses this issue by replacing the full affine-scaling
matrix with a diagonal approximation constructed from conjugacy-type conditions.
The resulting scaling strategy does not depend on gradient or subgradient
approximations, finite differences, or interpolation models. Instead, it uses
consecutive normalized recombination displacements in a conservative diagonal
update. Thus, the method does not assume that the generated directions are
noise-free; rather, it incorporates noisy directional information through a
bounded and regularized coordinatewise mechanism. Consequently, the method is
cheaper than full matrix-adaptation approaches and limited-memory affine-scaling
variants, while retaining a stable scaling mechanism for noisy large-scale
optimization problems.

The numerical experiments support the proposed design. In particular, the method
is robust under noisy function evaluations. The results suggest that conjugacy-based diagonal scaling provides
an effective compromise between stability, scalability, and computational cost. Future work may investigate noise-resistant definitions of the generated displacement direction \(\bd_k\). In particular, replacing the recombination-based direction in \eqref{dkMAES} by directions obtained from repeated evaluations, rank aggregation, trimmed recombination, or statistically reliable acceptance tests may further reduce the influence of ranking errors on the diagonal scaling update.

{\bf Supplementary Information}  The online version contains supplementary material available in \cite{suppMat}.

{\bf Funding} Morteza Kimiaei acknowledges financial support of the Austrian Science Foundation under \url{https://doi.org/10.55776/PAT2747625}.

\end{document}